\documentclass[11pt]{article}
\usepackage[utf8]{inputenc}
\usepackage[a4paper]{geometry}
\usepackage[english]{babel}
\usepackage[super]{nth}
\usepackage{amsmath, amsfonts, amsthm, physics, enumitem}
\usepackage{titlesec}
\usepackage{tikz} 
\usetikzlibrary{positioning}
\usepackage{mathtools}

\usetikzlibrary{patterns,arrows,decorations.markings}
\usepackage{authblk}
\usepackage{tikz-cd}
\usepackage{centernot}
\usepackage{caption}
\usepackage{float}

\newcommand{\KT}{\mathrm{KT}^{4}}

\newtheorem{thm}{Theorem}[section]
\newtheorem{prop}[thm]{Proposition}

\newtheorem{q}[thm]{Question}
\newtheorem{rmk}[thm]{Remark}

\pgfdeclareimage[height=15pt]{g_tick}{green_tick}
\pgfdeclareimage[height=15pt]{y_tick}{yellow_tick}
\pgfdeclareimage[height=15pt]{qu}{question}

\title{Almost K\"ahler Kodaira-Spencer problem}
\author{Tom Holt\thanks{Thomas.Holt@warwick.ac.uk}\,   and Weiyi Zhang\thanks{Weiyi.Zhang@warwick.ac.uk}}
\affil{Mathematics Institute,  University of Warwick, Coventry, CV4 7AL, England}

\date{}
\begin{document}

\maketitle

\begin{abstract}
We show that the almost complex Hodge number $h^{0,1}$ varies with different choices of almost K\"ahler metrics. This answers the almost K\"ahler version of a question of Kodaira and Spencer. 
\end{abstract}

\section{Introduction}
As motivated by Hodge theory for complex manifolds, it is important to study the space of $\bar\partial$-harmonic forms with any given Hermitian metric. Precisely, the almost complex structure $J$ on $M$ induces a decomposition of the complexified cotangent bundle $T^*M\otimes \mathbb C=(T^*M)^{1, 0}\oplus (T^*M)^{0,1}$, which in turn induces a decomposition of complex differential forms into $(p,q)$-forms. We define $\bar\partial$ (respectively $\partial$) to be the component of the exterior derivative $d$ that takes $(p,q)$-forms to forms of type $(p, q+1)$ (respectively $(p+1,q)$). Notice we no longer have $d=\partial+\bar\partial$ or $\bar\partial^2=0$ when $J$ is not integrable. Given an Hermitian metric we also define the operator $\bar{\partial}^{*}=-*\partial*$ along with the $\bar{\partial}$-Laplacian
$$\Delta_{\bar\partial}=\bar\partial\bar\partial^*+\bar\partial^*\bar\partial.$$
Here $*$ denotes the Hodge star.
The space $\mathcal H^{p,q}$ is defined to be the kernel of $\Delta_{\bar\partial}$ in the space of $(p, q)$-forms. When the manifold is compact, as the notation would suggest, $\bar{\partial}^{*}$ is the $L^2$ adjoint of $\bar\partial$ with respect to an Hermitian metric and so we have
  $\ker\Delta_{\bar\partial}=\ker\bar\partial\cap\ker\bar\partial^*$. Using the initial definition of $\bar{\partial}^{*}$, this is equivalent to $\mathcal H^{p,q}=\ker\bar\partial\cap \ker\partial *$. Since $\Delta_{\bar\partial}$ is an elliptic operator, the Hodge number $h^{p,q}=\dim\mathcal H^{p,q}$ is finite. 
  In this paper, we assume the manifold to be compact. 

For complex structures, $\mathcal H^{p,q}$ does not depend on the choice of the Hermitian metric since it is isomorphic to the Dolbeault cohomology group. In the almost complex setting, Kodaira and Spencer asked the following question which appeared as Problem 20 in Hirzebruch's 1954 problem list \cite{Hir}. For history and related problems, please see our paper \cite{HZ}.

\begin{q}[Kodaira-Spencer]\label{KS}
Let $M$ be an almost complex manifold. For any given Hermitian structure we can consider the numbers $h^{p,q}$. Are these numbers independent of the choice of the Hermitian structure?
\end{q}

When the answer is affirmative for some $(p,q)$, $h^{p,q}$ defines an almost complex invariant. Apparently, $h^{0,0}=h^{n,n}=1$ for any Hermitian structure on a connected $2n$-dimensional manifold $M$. It is true that $h^{p,0}$ are almost complex invariants \cite{CZ}. In fact, they are almost complex birational invariants shown for closed almost complex $4$-manifolds \cite{CZ2}. Moreover, as Serre duality  holds for $\mathcal H^{p,q}$ \cite{CZ}, we know $h^{p,n}=h^{n-p,0}$ are also almost complex invariants.

However, not all $h^{p,q}$ are almost complex invariants. Namely, Question \ref{KS} was answered negatively for $(p,q)=(0,1)$ on a $4$-manifold, the Kodaira-Thurston manifold, in \cite{HZ}. For all the examples constructed there, $h^{0,1}$ was only shown to differ in value when one of the two metrics used is almost K\"ahler and the other is not. It is then natural to ask for the answer to a version of Question \ref{KS} where the choice of Hermitian structure is required to be almost K\"ahler. In dimension $4$, although it is not known whether $h^{1,1}$ is an almost complex invariant, it is actually an almost K\"ahler invariant, {\it i.e.} independent of the choice of almost K\"ahler metrics compatible with $J$. More precisely, $h^{1,1}=b^-+1$ for almost K\"ahler metrics (Proposition 6.1 in \cite{HZ}). 

We summarise the known answers to the Kodaira-Spencer question for a closed 4-manifold in the following Hodge diamond picture.

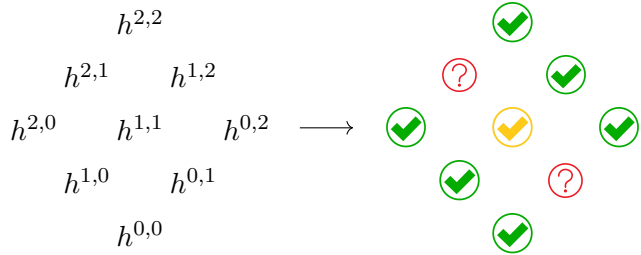
\begin{center}
\begin{figure}[H]
\centering
\captionsetup{justification=centering}
    \tikz [scale = .7,decoration={
       markings,
       mark=at position 1 with {\arrow[scale=3,black]{latex}};}] {
    
     \node (n1) at (0,0)  {$h^{0,0}$} ;
     \node (n2) at (1,1) {$h^{0,1}$};
     \node (n3) at (-1,1)  {$h^{1,0}$};
     \node (n4) at (0,2)  {$h^{1,1}$};
     \node (n5) at (2,2)  {$h^{0,2}$} ;
     \node (n6) at (-2,2) {$h^{2,0}$};
     \node (n7) at (-1,3)  {$h^{2,1}$};
     \node (n8) at (1,3)  {$h^{1,2}$};
     \node (n8) at (0,4)  {$h^{2,2}$};

    \draw [->] (3,2) -- (4,2);

    \node (n1) at (7,0)  {\pgfuseimage{g_tick}} ;
    \node (n2) at (8,1) {\pgfuseimage{qu}};
    \node (n3) at (6,1)  {\pgfuseimage{g_tick}};
    \node (n4) at (7,2)  {\pgfuseimage{y_tick}};
    \node (n5) at (9,2)  {\pgfuseimage{g_tick}} ;
    \node (n6) at (5,2) {\pgfuseimage{g_tick}};
    \node (n7) at (6,3)  {\pgfuseimage{qu}};
    \node (n8) at (8,3)  {\pgfuseimage{g_tick}};
    \node (n8) at (7,4)  {\pgfuseimage{g_tick}};
}
\caption*{  Green tick: Almost Hermitian metric invariant \\ Yellow tick: Almost K\"ahler metric invariant \\ Question mark: Not almost Hermitian metric invariant,  is it almost K\"ahler invariant?}
\end{figure}
\end{center}
Hence, in dimension $4$, the almost K\"ahler version of the Kodaira-Spencer question is reduced to asking whether $h^{0,1}=h^{2,1}$ is independent of the choice of almost K\"ahler metric.

\begin{q}[Question 5.2 in \cite{HZ}]\label{akKS}
Can we construct an almost complex structure $J$ on the Kodaira-Thurston, or more generally on a $4$-manifold, such that $h^{0,1}$ varies with different choices of almost K\"ahler metrics?
\end{q}

Since $h^{p,q}$ has been shown to be almost K\"ahler metric invariant for all except the one remaining case of $h^{0,1}$, we might expect to have a negative answer to Question \ref{akKS}. However, in this paper we answer this question affirmatively. 

\begin{thm}\label{main}
There exist almost complex structures on the Kodaira-Thurston manifold such that $h^{0,1}$ varies with different choices of almost K\"ahler metrics.
\end{thm}

In fact, the almost complex structures in the theorem are the same as the ones used in \cite{HZ}, which will be recalled in Section \ref{KTdef}. Although choosing a different family of metrics, we obtain a similar elliptic system, which could also be solved using our PDE-ODE-NT method introduced in \cite{HZ}. This method uses Fourier theory for the Heisenberg group to transform the elliptic PDE system of functions on the Kodaira-Thurston manifold to a set of countably many first order linear ODE systems and a set of countably many zeroth order linear equations. Note that these ODEs are on one variable functions by the Weil-Brezin transform. For convenience of the reader, we  develop the details of this harmonic analysis on the Kodaira-Thurston manifold in Section \ref{HAKT}. 
Then the ODE systems could be solved by analysing its Stokes phenomenon, while the zeroth order linear equations are reduced to the counting of lattice points on certain circles. In both cases, we can apply the results in \cite{HZ}.

Finally, we would like to comment a few more words on computing $h^{1,1}$ and determining whether it is an almost complex invariant by using our PDE-ODE-NT method. As shown by an example in \cite{HZ} as well as a recent computation of Tardini and Tomassini, the crux is to analyse the solutions of ODE systems. For the $h^{0,1}$ computations in both of our papers, we eventually reduce these ODE systems to confluent hypergeometric equations. However, we get $4\times 4$ systems for $h^{1,1}$, and a complete picture of their Stokes phenomenon seems much more difficult to achieve.

\section{Preliminary definitions}\label{KTdef}

The Kodaira-Thurston manifold $\KT$ is defined to be $\Gamma \backslash G$, where the group $G$ is given by $\mathbb{R}^4$ along with the group operation
\begin{equation}
\begin{pmatrix}
    t_0 \\
    x_0 \\
    y_0 \\
    z_0
\end{pmatrix} \circ
\begin{pmatrix}
    t\\
    x\\
    y\\
    z
\end{pmatrix}=
\begin{pmatrix}
    t+t_0\\
    x+x_0\\
    y+y_0\\
    z+z_0+x_0y
\end{pmatrix} \label{quot}\end{equation}
and $\Gamma = \mathbb{Z}^4$ is a discrete subgroup acting on $G$ from the left. It should be noted that $G$ is often given as $\mathbb{R} \times Nil^3 $ where 
$$Nil^3 = \left\{ \begin{pmatrix}
    1   &   x   &   z   \\
    0   &   1   &   y   \\
    0   &   0   &   1
\end{pmatrix}\middle| x,y,z\in \mathbb{R}\right\}.$$
We define the following left-invariant frame on $G$ 
$$e_1 = \partial_t \quad \quad e_2 = \partial_x \quad \quad e_3 = \partial_y + x \partial_z \quad \quad e_4 = \partial_z $$
along with its dual frame
$$ e^1 = dt \quad \quad e^2 = dx \quad \quad e^3 = dy \quad \quad e^4 = dz - x dy.$$
The structure equations for the dual frame are
$$de^1 = de^2 = de^3 = 0 \quad \quad de^4 = -e^2 \wedge e^3. $$
Note that although we defined this frame on $G$, since it is left-invariant it also induces a frame on $\KT$. 

We shall now consider a family of almost complex structures on $\KT$  given by the matrix
$$J_{a,b} = \begin{pmatrix}
        0 & -1 & 0 &  0\\
        1 &  0 & 0 &  0\\
        0 &  0 & a &  b\\
        0 &  0 & c &  -a
        \end{pmatrix}$$
acting on our frame, with $c = -\frac{a^{2}+1}{b}, a, b\in \mathbb R$.
We can then define the vector fields
$$V_{1} = \frac 12 \left(e_1 - ie_2 \right)\quad \mathrm{\&}\quad V_{2} = \frac 12 \left( e_3 -\frac{a-i}{b}e_4 \right) $$
spanning $T^{1,0}_{p}M$ at every point, along with their dual 1-forms
$$\phi^1 = e^1 +ie^2 \quad\mathrm{\&}\quad \phi^2 = (1-ai) e^3 -ib e^4. $$
Here the structure equations become $d\phi^1 = 0$ and
\begin{align*}
    d\phi^2 &= ib e^2 \wedge e^3\\
    &= ib \left(-\frac i2 (\phi^1 - \bar{\phi}^1) \right)\wedge \left( \frac 12 (\phi^2 + \bar{\phi}^2)\right)\\
    &= \frac b4 \left(\phi^{12} + \phi^{1\bar{2}}+ \phi^{2\bar{1}} - \phi^{\bar{1}\bar{2}} \right).
\end{align*}

We can obtain an almost K\"ahler structure if we also consider a family of Hermitian metrics $h_{a,b,\rho}$ compatible with $J_{a,b}$, along with its associated 2-form $\omega$ and complexified Riemannian metric $g$
$$h_{a,b,\rho} = 2(\phi^1 \otimes \bar{\phi}^1 + \rho \phi^2 \otimes \bar{\phi}^2 ), $$
$$g = \frac 12 ( h + \bar h ) =  \phi^1 \otimes \bar{\phi}^1 + \rho \phi^2 \otimes \bar{\phi}^2 + \bar{\phi}^1 \otimes \phi^1 + \rho \bar{\phi}^2 \otimes \phi^2,$$
\begin{align*}
    \omega = -\frac{i}{2}( h - \bar h ) &= -2i(\phi^1 \wedge \bar{\phi}^1 + \rho\phi^2 \wedge \bar{\phi}^2 )\\
    &= 4(e^2 \wedge e^1 + \rho e^3 \wedge e^4).
\end{align*}
This does indeed define an almost K\"ahler structure as 
\begin{align*}
    d\omega &= -4\rho e^3 \wedge de^4 =0.
\end{align*}
Since $V_1, \bar{V}_1 , \frac{1}{\sqrt \rho} V_2, \frac{1}{\sqrt \rho} \bar{V}_2$ are orthonormal with respect to $g$, we can define the volume form as $vol = \rho \phi^{12\bar{1}\bar{2}}$.

\section{Constructing equations}
Given the K\"ahler structure defined above, and for some fixed $a,b,\rho \in \mathbb{R}$, we write a general smooth $(0,1)$-form as $s = f\bar{\phi}^1 + g\bar{\phi}^2$ with $f,g\in C^{\infty}(\KT)$. 
We want to find which of these forms are $\bar{\partial}$-harmonic, this turns out to be identical to asking which $s$ satisfy both $\bar{\partial}s = 0$ and $\partial* s =0$. Here $*$ denotes the Hodge star operator defined by
$$\alpha \wedge * \bar{\beta} = g(\alpha,\beta) vol $$
for any two differential forms $\alpha, \beta$. In particular we have
$$*\bar{\phi}^1 = \rho \phi^{2\bar{1}\bar{2}}, \quad \quad *\bar{\phi}^2 = - \phi^{1\bar{1}\bar{2}}. $$

From the first of these conditions on $s$ we get
\begin{align*}
    \bar{\partial} (f\bar{\phi}^1 + g\bar{\phi}^2) &= -\bar{V}_2 (f) \phi^{\bar{1}\bar{2}}+ \bar{V}_1 (g) \phi^{\bar{1}\bar{2}} + g\bar{\partial}\bar{\phi}^2 \\
    &= \left(-\bar{V}_2 (f) + \bar{V}_1 (g) +g \frac b4 \right)\phi^{\bar{1}\bar{2}} =0,
\end{align*}
and from the second we get
\begin{align*}
    \partial * (f\bar{\phi}^1 + g\bar{\phi}^2) &= \partial (f\rho \phi^{2\bar{1}\bar{2}} - g\phi^{1\bar{1}\bar{2}} )\\
    &= \rho V_1 (f) \phi^{12\bar{1}\bar{2}} + f \rho \partial\phi^2 \wedge \phi^{\bar 1 \bar 2 } + f \rho  \phi^{2 \bar 1} \wedge \partial \phi^{\bar 2 } +V_2 (g)  \phi^{12\bar{1}\bar{2}} -g \phi^{1 \bar 1 } \wedge \partial \phi^{\bar 2 }\\
    &= \left(\rho V_1 (f) +f \frac {b\rho}{4} - f \frac {b\rho}{4} +V_2 (g)\right) \phi^{12\bar{1}\bar{2}}\\
    &= \left(\rho V_1 (f) +V_2 (g)\right) \phi^{12\bar{1}\bar{2}} =0.
\end{align*}
So we see that $f$ and $g$ must satisfy the following system of equations:
\begin{align}
-\bar{V}_2 (f) + \bar{V}_1 (g) +g \frac b4 &=0, \label{V1}\\
\rho V_1 (f) +V_2 (g) &=0. \label{V2}
\end{align}

\section{Harmonic Analysis on $\KT$}\label{HAKT}

In order to solve the system of equations \eqref{V1} and \eqref{V2} we now introduce a decomposition of $L^2(\KT)$ into the direct sum of simpler classes of functions, denoted $\mathcal{H}_{k,m,n}$ and $\mathcal{H}_{k,l,m,0}$. We will then see that all solutions of our system of equations can be found by looking for solutions within $\mathcal{H}_{k,m,n}\cap C^{\infty}(\KT)$ and $\mathcal{H}_{k,l,m,0}\cap C^{\infty}(\KT)$.

\begin{prop}
The space of square-integrable functions on the Kodaira-Thurston manifold $\KT$ decomposes in the following way:
$$L^2 (\KT) = \widehat{\bigoplus_{k\in \mathbb{Z}}}\left(\widehat{\bigoplus}_{l,m \in \mathbb{Z}}\mathcal{H}_{k,l,m,0}\right) \oplus \left( \widehat{\bigoplus}_{\substack{n\in \mathbb{Z}\backslash \{0\}, k\in \mathbb Z \\ m\in \{0,\,1,\,\dots ,\, \abs{n}-1\}}}\mathcal{H}_{k,m,n} \right)$$
where the symbol $\widehat{\bigoplus}$ denotes the closure with respect to the $L^2$ norm of the direct sum and we define
$$\mathcal{H}_{k,l,m,0}=\left\{f_{k,l,m,0}\, e^{2\pi i (kt+lx+my)}\in L^{2}(\KT) \, \middle| \, f_{k,l,m,0}\in \mathbb{C} \right\}, $$
$$\mathcal{H}_{k,m,n}=\left\{ \sum_{\xi \in \mathbb{Z}}f_{k,m,n}(x+\xi) e^{2 \pi i (kt+(m+n\xi)y +nz)} \in L^2 (\KT) \, \middle| \, f_{k,m,n}(x) \in L^2 (\mathbb{R}) \right\}. $$
    \begin{proof}
    Since $L^2(\KT)$ is the closure of $C^{\infty}(\KT)$ with respect to the $L^2$ norm, this means we only need to show that all \textit{smooth} functions can be written as a sum of functions contained in $\mathcal{H}_{k,m,n}$ and $\mathcal{H}_{k,l,m,0}$. If this is the case then $C^{\infty}(\KT)$ must be contained in the closure of the direct sum of these spaces, and thus their closure must consist of the whole of $L^2(\KT)$.  
    
    Recall the definition of $\KT$ as the group $G=(\mathbb{R}^4,\circ)$ modulo the left action of $\Gamma =\mathbb{Z}^4$. In this way we can consider functions on $\KT$ as functions on $\mathbb{R}^4$ satisfying the condition
    \begin{equation}
        f(t,x,y,z) = f(t+\alpha, x+\xi, y+\beta, z+\gamma+\xi y)  \label{condition}
    \end{equation}
    for all $\alpha, \beta, \gamma, \xi \in \mathbb{Z}^4$. In particular, when $\xi=0$ this condition tells us the functions must be periodic with respect to $t, y$ and $z$.
    From classical Fourier analysis we can therefore write any smooth function $f\in C^{\infty}(\mathbb{R}^4)$ satisfying \eqref{condition} as the series
    $$f(t,x,y,z) = \sum_{k,m,n\in \mathbb{Z}}f_{k,m,n}(x)e^{2 \pi i (kt+my+nz)} $$
    where
    $$f_{k,m,n}(x) = \int_{[0,1]^3}f(t,x,y,z)e^{-2 \pi i (kt+my+nz)}dt dy dz. $$
    This almost gives us a decomposition of functions on $\KT$ but the terms in the above series may not individually satisfy \eqref{condition} and thus may not be well defined functions on $\KT$. To fix this we need to take into account the case in \eqref{condition} when $x_0$ is allowed to be non-zero, that is to say we need 
    $$f(t,x,y,z) = f(t,x+\xi , y , z+\xi  y)$$
    for all $\xi \in \mathbb{Z}$. After taking the Fourier expansion of both sides this becomes
    $$\sum_{k,m,n\in \mathbb{Z}}f_{k,m,n}(x)e^{2 \pi i (kt+my+nz)} = \sum_{k,m,n\in \mathbb{Z}}f_{k,m,n}(x+\xi)e^{2 \pi i (kt+my+n(z+\xi y))}, $$
    or equivalently
    $$\sum_{k,m,n\in \mathbb{Z}}f_{k,m,n}(x)e^{2 \pi i (kt+my+nz)} = \sum_{k,m,n\in \mathbb{Z}}f_{k,m-n \xi,n}(x+\xi)e^{2 \pi i (kt+my+nz)}. $$
    Using the fact that Fourier coefficients of smooth functions are unique, this tells us that we need $f_{k,m+n\xi,n}(x) = f_{k,m,n}(x+\xi)$. In fact, since we are looking at smooth functions on $\KT$, not just continuous functions, we require that 
    $$\partial_{x}^{i}f(t,x,y,z) = \partial_{x}^{i}f(t,x+\xi,y,z+\xi y)$$
    and therefore $\frac{d^i}{dx^i}f_{k,m+n\xi,n}(x) = \frac{d^i}{dx^i}f(x+\xi)$ for all $i \in \mathbb{N}$ and all $x_0 \in \mathbb{Z}$.
    
    We now find that we have two cases. When $n=0$ the above condition is just the requirement that $f_{k,m,0}(x)$ be periodic with respect to $x$ with period $1$, which means we can write
    $$f_{k,m,0}(x) = \sum_{l\in \mathbb{Z}}f_{k,l,m,0}e^{2\pi i lx} $$
    where
    $$f_{k,l,m,0} = \int_{[0,1]^4}f(t,x,y,z)e^{-2\pi i (kt+lx+my+nz)}dt dx dy dz. $$
    If instead $n\neq 0$ we have collections of functions $\{f_{k,m+n\xi,n}(x) \}_{\xi \in \mathbb{Z}}$ which are just translations of each other. This means for any $k\in\mathbb{Z}$ we can write
\begin{align*}
    \sum_{\substack{m,n\in \mathbb{Z}\\ n\neq 0}}f_{k,m,n}(x)e^{2\pi i (kt+my+nz)} &= \sum_{\substack{n \in \mathbb{Z}\backslash \{0\} \\ m\in \{0,1,\dots , \, \abs{n}-1\}}} \sum_{\xi \in \mathbb{Z}}f_{k,m+n\xi,n}(x)e^{2\pi i (kt+(m+n\xi)y + nz)}\\
    &= \sum_{\substack{n \in \mathbb{Z}\backslash \{0\}\\ m\in \{0,1,\dots , \, \abs{n}-1\}}} \sum_{\xi \in \mathbb{Z}}f_{k,m,n}(x+\xi)e^{2\pi i (kt+(m+n\xi)y + nz)}.
\end{align*}
    Notice that 
    $$f=\sum_{\xi \in \mathbb{Z}}f_{k,m,n}(x+\xi)e^{2\pi i (kt+(m+n\xi)y + nz)}$$ 
    does indeed satisfy $\eqref{condition}$ and thus is well defined on $\KT$. For any $k,m,n\in \mathbb{Z}$ with $n\neq 0$ we have the map 
    $W_{k,m,n}:L^2(\mathbb{R}) \rightarrow \mathcal{H}_{k,m,n}$ defined by  $$W_{k,m,n}:f_{k,m,n}(x)\mapsto \sum_{\xi \in \mathbb{Z}}f_{k,m,n}(x+\xi)e^{2\pi i (kt+(m+n\xi)y + nz)}. $$
    This is essentially the Weil-Brezin transform \cite{Aus, Fol}, mapping $\mathcal{S}(\mathbb{R})$ bijectively onto $\mathcal{H}_{k,m,n}\cap C^{\infty}(\KT)$ where $\mathcal{S}(\mathbb{R})$ denotes the space of Schwartz functions.
    $$\mathcal{S}(\mathbb{R}) = \left\{f\in C^{\infty}(\mathbb{R})\, : \, \sup_{x}\abs{x^p \frac{d^s}{dx^s}f(x)}<\infty \quad  \forall p,s\in \mathbb{N} \right\}. $$
    We can see that this is the case by noting that if we have a function $f\in \mathcal{H}_{k,m,n}\cap C^{\infty}(\KT)$, then all its derivatives are at least locally bounded. In particular, their Fourier coefficients are locally bounded, meaning any $x_0 \in \mathbb{R}$ has a neighbourhood $N_{x_0}$ such that
    \begin{equation*}
        \sup_{\substack{k,m,n\in \mathbb{Z}\\x\in N_{x_0}}}\abs{k^p m^q n^r \frac{d^s}{dx^s} f_{k,m,n}(x)}< \infty 
    \end{equation*}
    for all $p,q,r,s \in \mathbb{N}$ and so $f_{k,m,n}(x)$ must be Schwartz.
    Conversely if a collection of Fourier coefficients satisfies the above condition, then the series they define converges. Furthermore, all the derivatives of the terms of the Fourier series converge to zero locally uniformly and so the series is infinitely differentiable. 
    As a result $\mathcal{H}_{k,m,n}\cap C^{\infty}(\KT)$ is exactly the image of $\mathcal{S}(\mathbb{R})$ under $W_{k,m,n}$.
    
    Combining the two cases described above gives us a full decomposition for smooth functions on $\KT$
    \begin{align*}
        f = \sum_{k\in \mathbb{Z}}\bigg( &\sum_{l,m\in \mathbb{Z}}f_{k,l,m,0}e^{2\pi i (kt + lx + my)}\\
        + &\sum_{\substack{n\in \mathbb{Z}\backslash \{0\}\\ m\in \{0,1,\dots, \, \abs{n}-1\}}}\sum_{\xi\in \mathbb{Z}} f_{k,m,n}(x+\xi)e^{2\pi i (kt + (m+n\xi)y + nz)} \bigg)
    \end{align*}    
    with $f_{k,l,m,0} \in \mathbb{C}$ and $f_{k,m,n}(x) \in \mathbb{S}(\mathbb{R})\subset L^2(\mathbb{R})$. 
    
    \end{proof}
\end{prop}

\begin{rmk}
From a representation theory standpoint this decomposition of $L^2(\KT)$ corresponds to the decomposition of the regular representation $R$ of $G=(\mathbb{R}^4 , \circ)$ on $L^2(\KT)$
$$R(t_0,x_0,y_0,z_0)f(t,x,y,z) = f((t,x,y,z)\circ(t_0,x_0,y_0,z_0)) $$
in terms of irreducible representations
$$\sigma_{k,l,m}(t_0,x_0,y_0,z_0)= e^{2\pi i (kt_0+lx_0+my_0)} $$
acting on $\mathbb{C}$ with $k,l,m\in \mathbb{Z}$ and
$$\rho_{k,n}(t_0,x_0,y_0,z_0)f(x) = e^{2\pi i (kt_0 + n(z_0 + x y_0  )}f(x+x_0) $$
acting on $L^2(\mathbb{R})$ with $n\in \mathbb{Z}\backslash \{0\}$. Here $\sigma_{k,l,m}$ corresponds to the $n=0$ case and $\rho_{n}$ corresponds to the $n\neq 0$ case. 

$\KT$ may be written as the direct product of $S^1$ with the Heisenberg manifold the representations, and so $\sigma_{k,l,m}$ and $\rho_{k,n}$ can be derived from similar representations on the Heisenberg manifold which are given by virtue of the Stone-von Neumann theorem. For more detail see \cite{Aus, Fol}.
\end{rmk}

It is actually possible to generalise the above proposition to apply to any $d$-dimensional torus bundle over $S^1$. We can define such a bundle by taking the trivial bundle $\mathbb{T}^d \times[0,1]$ and identifying the ends using a torus automorphism given by $A \in GL_d(\mathbb{Z})$, for instance in the case of $\KT$ we have $A=\begin{pmatrix}
    1   &   0   &   0\\
    0   &   1   &   1\\
    0   &   0   &   1
\end{pmatrix}$.
The two cases $n=0$ and $n\neq 0$ from the Proposition turn out to correspond to the cases when the orbit of the group generated by $A$, acting on an element of $\mathbb{Z}^d$, has either finitely many or infinitely many elements.

\section{Solving the system of PDEs}

We can now return to solving our system of PDEs, equipped with the technique we introduced in the previous section. It is easy to check that $e_1, e_2, e_3$ and $e_4$, as defined earlier, all preserve the spaces $\mathcal{H}_{k,l,m,0}\cap C^{\infty}(\KT) = \mathcal{H}_{k,l,m,0} $ and $\mathcal{H}_{k,m,n} \cap C^{\infty}(\KT)$. This means that for any differential operator $V$ given by a linear combination of $e_1, e_2, e_3$ and $e_4$, and any $f\in C^{\infty}(\KT)$, if $f$ decomposes into the sum of $F_{k,l,m,0}\in \mathcal{H}_{k,l,m,0}$ and $F_{k,m,n}\in \mathcal{H}_{k,m,n}\cap C^{\infty}(\KT)$ then $V(f)$ decomposes into the sum of $V (F_{k,l,m,0})\in \mathcal{H}_{k,l,m,0}$ and $V(F_{k,m,n})\in \mathcal{H}_{k,m,n}\cap C^{\infty}(\KT)$. Since $V_1$ and $V_2$ are two such differential operators we can similarly decompose the equations \eqref{V1} and \eqref{V2}. 
By the uniqueness of terms in the decomposition,  every term on the left hand side of both equations is zero. 

We therefore see that the $n=0$ case gives us independent solutions 
$$ f = f_{k,l,m,0}\, e^{2\pi i (kt+lx+my)}\in \mathcal{H}_{k,l,m,0}$$ 
$$g = g_{k,l,m,0}\, e^{2\pi i (kt+lx+my)}\in \mathcal{H}_{k,l,m,0}$$
for every $k,l,m\in \mathbb{Z}$ such that $f_{k,l,m,0},g_{k,l,m,0}\in \mathbb{C}$ solve
$$-m f_{k,l,m,0} + \left(k+il -i\frac{b}{4\pi} \right)g_{k,m,n}=0, $$
$$\rho \left(k-il \right)f_{k,m,n} + mg_{k,m,n} =0.$$
We also see that the $n\neq 0$ case gives us independent solutions
$$f= \sum_{\xi \in \mathbb{Z}}f_{k,m,n}(x+\xi) e^{2 \pi i (kt+(m+n\xi)y +nz)} \in \mathcal{H}_{k,m,n}\cap C^{\infty}(\KT) $$
$$g= \sum_{\xi \in \mathbb{Z}}g_{k,m,n}(x+\xi) e^{2 \pi i (kt+(m+n\xi)y +nz)} \in \mathcal{H}_{k,m,n}\cap C^{\infty}(\KT) $$
for every $k,m,n\in \mathbb{Z},\, n\neq 0$ such that $f_{k,m,n}(t), g_{k,m,n}(t)\in \mathcal{S}(\mathbb{R})$ solve
$$-\left(m+xn - \frac{a+i}{b}n\right)f_{k,m,n} + \left(k+\frac{1}{2\pi}\frac{d}{dx} -i\frac{b}{4\pi} \right)g_{k,m,n}=0, $$
$$\rho \left(k-\frac{1}{2\pi}\frac{d}{dx} \right)f_{k,m,n} + \left( m+nx - \frac{a-i}{b}n\right)g_{k,m,n} =0.$$

In this way we can split the hodge number $h^{0,1}$ into two parts
$$h^{0,1} = h_{0,1}'+h_{0,1}'' $$
where $h_{0,1}$ counts the number of solutions arising from the $n=0$ case and $h_{0,1}''$ counts the number of solutions arising from the $n\neq 0$ case.

\subsection{Solving the $n=0$ case}\label{n0}
Here we are looking for solutions $f_{k,l,m,0},g_{k,l,m,0}\in \mathbb{C}$ to the system of Diophantine equations
$$-m f_{k,l,m,0} + \left(k+il -i\frac{b}{4\pi} \right)g_{k,m,n}=0 $$
$$\rho \left(k-il \right)f_{k,m,n} + mg_{k,m,n} =0$$
When $m=0$ we can see that either $k=l=0$ and $g_{k,l,0,0}=0$ or $k=0, l=\frac{b}{4\pi}$ and $f_{k,l,0,0}=0$. These two cases give us the solutions
$$f=C_{1}, g=0 \quad \text{and} \quad f=0, g=C_{2}e^{2\pi i \frac{b}{4\pi}x} $$
for any $C_{1}, C_{2} \in \mathbb{C}$. 
Note that the second of these solutions is only possible when $\frac{b}{4\pi} \in \mathbb{Z}$.
Now setting $m\neq 0$ we can rearrange our equations to get
$$g_{k,l,m,0}=-\rho \frac{k-il}{m}f_{k,l,m,0}, $$
$$\left(m^2 + \rho k^2 + \rho l^2 - \frac{b\rho}{4\pi }(l+ki)\right)f_{k,l,m,0}=0. $$
Here we can see there will be solutions only when $k=0$ and non-zero $l,m\in \mathbb{Z}$ are chosen such that
$$m^2 + \rho l^2 -\frac{b}{4\pi}\rho l=0. $$
Such a choice for $l,m$ will give us the solutions
$$f= mC_{3}e^{2\pi i (lx+my)}, g=i\rho l C_{3}e^{2\pi i (lx+my)}  $$
$C_{3}\in\mathbb{C}$. But how many integer values for $l$ and $m$ satisfy this condition? Setting $d = \frac{b}{8\pi}$ we can restate the condition as
$$\left(\frac{m}{\sqrt{\rho}}\right)^2 +(l-d)^2 = d^2. $$
This is equivalent to asking how many points on the lattice given by $\mathbb{Z}\times\frac{1}{\sqrt \rho}\mathbb{Z}$ intersect a circle of radius $d$ and centre $(d,0)$. Note that the two solutions found when $m=0$ can be considered to correspond to the lattice points $(0,0)$ and $(\frac{b}{4\pi},0)$.

\subsection{Solving the $n\neq 0$ case}\label{no0}
Here we are looking for solutions $f_{k,m,n}(t), g_{k,m,n}(t)\in \mathcal{S}(\mathbb{R})$ to the system of ODEs
$$-\left(m+xn - \frac{a+i}{b}n\right)f_{k,m,n} + \left(k+\frac{1}{2\pi}\frac{d}{dx} -i\frac{b}{4\pi} \right)g_{k,m,n}=0, $$
$$\rho \left(k-\frac{1}{2\pi}\frac{d}{dx} \right)f_{k,m,n} + \left( m+nx - \frac{a-i}{b}n\right)g_{k,m,n} =0.$$
We can rearrange this into the form
$$\frac{d}{dx} \begin{pmatrix}
    f_{k,m,n}\\
    g_{k,m,n}
\end{pmatrix}=
(Ax+B)\begin{pmatrix}
    f_{k,m,n}\\
    g_{k,m,n}
\end{pmatrix}$$
where
$$A = 2\pi n\begin{pmatrix}
    0   &   \frac{1}{\rho}\\
    1   &   0
\end{pmatrix}, \quad \quad 
B = 2\pi\begin{pmatrix}
    k   &   \frac {1}{\rho}\left(m-\frac{a-i}{b}n\right)\\
    m - \frac{a+i}{b}n    &   i\frac{b}{4\pi}-k
\end{pmatrix}.$$

We now make use of the following theorem from \cite{HZ}.
\begin{thm}
Let $A,B \in M_{2}(\mathbb{C})$ be matrices and let $A$ have two distinct, real eigenvalues $\lambda_{1}$, $\lambda_{2}$ with $\lambda_{1}>0>\lambda_2$ then the equation
\begin{align}
    \frac{d}{dx}
    \begin{pmatrix}
        f\\
        g
    \end{pmatrix}
    = (Ax+B)
    \begin{pmatrix}
        f\\
        g
    \end{pmatrix}\label{Ax+B}
\end{align}
has a pair of solutions $f,g\in L^{2}(\mathbb{R})$ if and only if the following holds: Given $P\in GL(2, \mathbb C)$ such that $PAP^{-1}$ is diagonal and writing $PBP^{-1}$ as $\begin{pmatrix}
    b_{1} & b_{2}\\
    b_{3} & b_{4}
\end{pmatrix}$ we have $b_{2}b_{3}\in (\lambda_{1}-\lambda_{2})\cdot \mathbb{Z}^-$, and in this situation both $f$ and $g$ are Schwartz functions.
\end{thm}
In our case if we choose 
$$P= \frac{\sqrt 2}{2}\begin{pmatrix}
    \sqrt \rho  &   1\\
    \sqrt \rho  &   -1
\end{pmatrix} $$
then we find that 
$$PAP^{-1} = \frac{2\pi n}{\sqrt \rho}\begin{pmatrix}
    1   &   0\\
    0   &   -1
\end{pmatrix}, \quad 
PBP^{-1} = 2\pi \begin{pmatrix}
    \frac{1}{\sqrt \rho}\left(m -\frac{an}{b}\right)+\frac{b}{8\pi}i & k-\frac{n}{b\sqrt \rho}i-\frac{b}{8\pi}i\\
    k+\frac{n}{b\sqrt \rho}i-\frac{b}{8\pi}i & -\frac{1}{\sqrt \rho}\left(m- \frac{an}{b}\right)+\frac{b}{8\pi}i
\end{pmatrix}.$$
So we have solutions only when 
$$4\pi^2 \left(k-\frac{n}{b\sqrt \rho}i-\frac{b}{8\pi}i\right)\left( k+\frac{n}{b\sqrt \rho}i-\frac{b}{8\pi}i\right)\in \frac{4\pi n}{\sqrt \rho}\mathbb{Z}^- .$$
The imaginary part of the left hand side is $-kb\pi $ so we are forced to set $k=0$, this leaves us with the condition that for some $u\in \mathbb{Z}^-$ 
$$b^4 \rho + 64\pi n u b^2 \sqrt \rho - 64 n^2 \pi ^2 =0, $$
or alternatively if we set $d=\frac{b}{8\pi}$
$$64\pi^2 \rho d^4 +64 \pi nu\sqrt \rho d^2-n^2=0. $$
But, since $\pi$ is transcendental, there can be no integer choice of $n$ and $u$ satisfying this condition unless $d$ and $\rho$ are chosen such that $8\pi \sqrt \rho d^2$ is a quadratic integer. In particular, $\sqrt \rho d^2 \in \frac{1}{\pi}\bar{\mathbb{Q}}$, where $\bar{\mathbb{Q}}$ is used to denote the algebraic numbers.

\subsection{Proof of Theorem \ref{main}}

Now we can finish the proof of Theorem \ref{main}. Our computation is workable for any parameters $d$ and $\rho$ with the results in \cite{HZ}. But here we are satisfied to exploit computation for special values of $d$ and $\rho$ in order to show $h^{0,1}$ could vary. 

Consider the case when $\frac{b}{8\pi}=d=1$. Bringing together the results from Sections \ref{n0} and \ref{no0} we find that
$$h_{0,1}' = \begin{cases}
    4 & \sqrt \rho \in \mathbb{Z},\\
    2 & \sqrt \rho \not\in \mathbb{Z},
\end{cases}$$
while $h_{0,1}''=0$ unless $\sqrt \rho\in \frac{1}{\pi} \bar{\mathbb{Q}}$. This clearly indicates that $h^{0,1}=h_{0,1}'+h_{0,1}''$ does not take one fixed value as we have
$$ h^{0,1} = \begin{cases} 
      4 & \sqrt \rho \in \mathbb{Z} ,\\
      2 & \sqrt \rho \in \mathbb{Q}\backslash \mathbb{Z} .
   \end{cases} $$
  
 This finishes the proof of  Theorem \ref{main}, {\it i.e.} $h^{0,1}$ varies under almost K\"ahler metric deformations.

Finally, we remark that we can actually achieve arbitrarily large values of $h^{0,1}$ by varying the almost complex structure $J_{a,b}$, as a result of the computation in \cite{HZ}.

\end{document}